# Non-Gaussian dynamics of a tumor growth system with immunization *

Mengli Hao[1,2], Ting Gao[2,3], Jinqiao Duan[2,3] and Wei Xu[1]
1. Department of Applied Mathematics
Northwest Polytechnology University, Xi'an, 710129, China
E-mail: haomengli@gmail.com
2. Institute for Pure and Applied Mathematics
University of California Los Angeles, CA 90095, USA
E-mail: jduan@ipam.ucla.edu
3. Department of Applied Mathematics
Illinois Institute of Technology, Chicago, IL 60616, USA
E-mail: tinggao0716@gmail.com

June 2, 2012

**Abstract**

This paper is devoted to exploring the effects of non-Gaussian fluctuations on dynamical evolution of a tumor growth model with immunization, subject to non-Gaussian $\alpha$-stable type Lévy noise. The corresponding deterministic model has two meaningful states which represent the state of tumor extinction and the state of stable tumor, respectively. To characterize the lifetime for different initial densities of tumor cells staying in the domain between these two states and the likelihood of crossing this domain, the mean exit time and the escape probability are quantified by numerically solving differential-integral equations with appropriate exterior boundary conditions. The relationships between the dynamical properties and the noise parameters are examined. It is found that in the different stages of tumor, the noise parameters have different influence on the lifetime and the

*This work was done while Mengli Hao was visiting the Institute for Pure and Applied Mathematics (IPAM), Los Angeles, CA 90025, USA. This work was supported by the NSFC Grants 10971225, 10932009 and 11172233, the NSF Grants 0731201 and 1025422, and the Fundamental Research Funds for the Central Universities (HUST 2010ZD037).

likelihood inducing tumor extinction. These results are relevant for determining efficient therapeutic regimes to induce the extinction of tumor cells.

# 1 Introduction

Recently, with the increase in the number of people with tumor, the growth law of tumor cells has attracted a lot of attention. Although surgery and chemo- and radio-therapies are the main treatment strategies, they cannot always cure the cancer completely. Many patients who have shown significant improvement still have to face the risk of relapse. Therefore, more efficient therapies with fewer side effects are urgently called for. Since the immune reaction is of importance for producing clinically observed phenomena, such as tumor dormancy, oscillations in tumor size and spontaneous tumor regression, efforts have been devoted to studying the immunotherapy and the mechanisms of interplay between tumor tissue and immune cells [8,13,14,16]. A model with immunization that can represent interaction and competition between tumor tissue and immune system has been proposed based on the Michaelis-Menten mechanism, which is one of the most important chemical reaction mechanisms in biochemistry. The response of tumor tissue to immunization process can be expressed by means of a predator-prey model. Tumor cells and immune cells act as the roles of "preys" and "predators", respectively.

Due to the influence of time-variability of the environment, such as temperature variations or changes in the local concentrations of biochemical agents, the system is inevitably impacted by random fluctuations [9]. Thus, it is worthwhile to pursue the theoretical study of the roles of the external noise in this system to obtain insight into the reaction mechanisms.

In the existing work, certain dynamical properties of this model under Gaussian noise (in terms of Brownian motion) have been studied, such as stability, stochastic resonance, mean first passage time [1,4,9,10,12,17,20,21, 26–28]. These results have established that environment fluctuations affect the growth and extinction of the tumor cells and thus are crucial for further understanding the model.

In the present work, we are interested in the dynamical behaviors of the tumor cell growth model with immunization subject to the non-Gaussian noise (in terms of Lévy motion). We use the mean exit time and escape probability to quantify the evolution of tumor cell growth, including the examination of when the tumor cells may be eliminated with higher probability.

This paper is organized as follows. Section 2 presents the tumor growth model with immunization, under the influence of non-Gaussian Lévy noise. In section 3, we provide a brief introduction to Lévy processes. In section 4, we recall the mathematical formulation for the mean exit time and escape probability, together with the numerical algorithms. Numerical results based on mean exit time and escape probability are presented in section 5. The paper is ended by a discussion in the final section.

## 2 Tumor growth model driven by Lévy noise

The tumor growth model in the presence of immunization can be expressed by a predator-prey model, in which the tumor cells act as preys, and the cytotoxic cells as predators. The process of the reactions between the tumor tissue and the immune cells can be well approximated by a saturating, enzymatic-like process whose time evolution equations follow the catalytic Michaelis-Menten kinetics [10]. So the formulation of the model is equivalent to the following reaction scheme:

$$X \xrightarrow{\lambda} 2X, \qquad (1a)$$

$$X + Y \xrightarrow{k_1} Z \xrightarrow{k_2} Y + P, \qquad (1b)$$

where $X, Y, Z$ and $P$ represent the populations of tumor cells, active cytotoxic cells, the complex of tumor cells and cytotoxic cells and dead or non-replicating tumor cells, respectively. The reaction process can be explained as follows: First, the tumor cells proliferate spontaneously at a rate $\lambda$ (Here, we ignore the transformation of normal cells into neoplastic ones because its rate is very low compared with the rate of neoplastic cell replication [13, 16]). Then, the cytotoxic cells bind the tumor cells to the complex $Z$ with the kinetic constant $k_1$. At last, the complex $Z$ dissociates at a rate proportional to $k_2$. In the above presentation, $Y$ behaves as the enzymes in the Michaelis-Mentecn reaction, so $Y + Z = E$ can be taken as a constant. Typical experimental values of parameters are: $\lambda$ =0.2-1.5 $day^{-1}$, $k_1$ =0.1-18 $day^{-1}$, $k_2$ =0.2-18 $day^{-1}$ [9, 16].

In the limit case, the production of X-type cells inhibited by a hyperbolic activation is the slowest process. By the quasi-steady-state approximation we can reduce the dimension of the chemical master equation and rewrite the resulting kinetics in the form of the first order differential equation [9, 16, 21]:

$$\frac{\mathrm{d}x}{\mathrm{d}t} = x(1 - \theta x) - \beta \frac{x}{x+1}, \tag{2}$$

with the potential function

$$U(x) = -\frac{x^2}{2} + \frac{\theta x^3}{3} + \beta x - \beta ln(x+1), \tag{3}$$

where $x$ is the normalized molecular density of tumor cells with respect to the maximum tissue capacity. And we use the following scaling formulas in the process of nondimensionalization:

$$x = \frac{k_1}{k_2} X, \theta = \frac{k_2}{k_1}, \beta = \frac{k_1 E}{\lambda}, t = \lambda t'.$$

The deterministic dynamical system Eq. (2), for some parameters, has two stable states and one unstable state (see [16, 21]). Namely, for $\theta < 1, 0 < \beta < \frac{(1+\theta)^2}{4\theta}$, the potential function $U(x)$ has two minima: $x_1$, $x_3$ and one maximum at $x_2$ (see Figure 1). That is to say, this system has three meaningful steady states:

$$\begin{aligned} x_1 &= 0, \\ x_2 &= \frac{1 - \theta - \sqrt{(1+\theta)^2 - 4\beta\theta}}{2\theta}, \\ x_3 &= \frac{1 - \theta + \sqrt{(1+\theta)^2 - 4\beta\theta}}{2\theta}. \end{aligned} \tag{4}$$

Without random fluctuations, system states finally approach one of the two stable states: (i) either the stable state $x_1 = 0$, where no tumor cells are present, namely, the tumor-free state (or the state of tumor extinction), (ii) or the other stable state $x_3$, where the tumor cell density does not increase but stays at a certain constant level, namely, the state of stable tumor.

However, the growth rate of tumor tissue is inevitably influenced by many environmental factors, such as the supply of nutrients, the immunological state of the host, chemical agents, temperature, and radiations [9, 16]. In the past, much work has been done on this model driven by Gaussian noise [1, 4, 9, 10, 12, 17, 20, 21, 26–28]. To describe the environment fluctuations more exactly, in this paper, we analyze the more realistic model including source of stochastic fluxes represented by non-Gaussian noise, modeled as time derivative of a Lévy motion $L_t$. Thus, we consider the model equation rewritten as follows:

$$\mathrm{d}X_t = f(X_t)\mathrm{d}t + \mathrm{d}L_t, \quad X(0) = x, \tag{5}$$

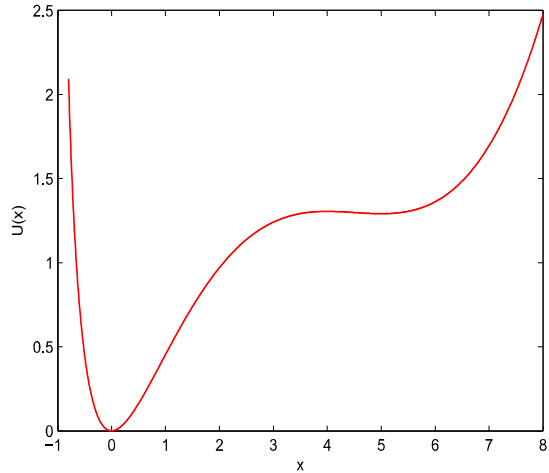

Figure 1: The potential function $U(x)$ for $\theta = 0.1, \beta = 3.0$.

where
$$f(X_t) = X_t(1 - \theta X_t) - \beta \frac{X_t}{X_t + 1},$$
and $L_t$ is a scalar Lévy motion defined in a probability space $(\Omega, \mathcal{F}, P)$. As will be discussed in the next section, the generating triplet for this Lévy motion is taken to be $(0, a, \varepsilon \nu_\alpha)$, i.e., zero drift coefficient, diffusion coefficient $a \geq 0$ and jump measure $\varepsilon \nu_\alpha(\mathrm{d}y)$, with $0 < \alpha < 2$. Here $\varepsilon \in [0, 1]$, being the noise intensity. According to the Lévy-Itô decomposition, this Lévy motion is the independent sum of a Brownian motion and a symmetric $\alpha$-stable process. The noise parameters $a, \alpha, \varepsilon$ means the diffusion coefficient or variation of the Gaussian noise part, stability index and the intensity of the non-Gaussian noise part, respectively.

Due to the addition of external noise, the number of tumor cells may fluctuate between the state of tumor extinction and the state of the stable tumor. Therefore, in this paper, we are interested in how the tumor density progresses in the range between the state of tumor extinction and the stable tumor state, denoted by $D = (x_1, x_3)$. We discuss the lifetime that the density of tumor cells stays in the range $D$ and the likelihood it gets out of the range to the left side, i.e., becoming extinction.
55

# 3 Lévy motion

Let $L = (L(t), t \geq 0)$ be a stochastic process defined on a probability space $(\Omega, \mathcal{F}, P)$. We say that it has independent increments if for each $n \in \mathbb{N}$ and each $0 \leq t_1 < t_2 < \cdots < t_{n+1} < \infty$ the random variables $(L(t_{j+1}) - L(t_j), 1 \leq j \leq n)$ are independent and that it has stationary increments if each $L(t_{j+1}) - L(t_j) \stackrel{d}{=} L(t_{j+1} - t_j) - L(0)$.

We say that $L(t)$ is a Lévy process if:
**(L1)** $L(0) = 0$ (a.s.);
**(L2)** $L(t)$ has independent and stationary increments;
**(L3)** $L(t)$ is stochastically continuous, i.e. for all $a > 0$ and for all $s > 0$

$$\lim_{t \to s} P(|L(t) - L(s)| > a) = 0.$$

Note that in the presence of (L1) and (L2), (L3) is equivalent to the condition
$$\lim_{t \to 0} P(|L(t)| > a) = 0$$
for all $a > 0$ (See [3]).

It is known that a scalar Lévy motion is completely determined by the Lévy-Khintchine formula (see [3, 23]).

**The Lévy-Khintchine Formula.** If $L = (L(t), t \geq 0)$ is a Lévy motion, then its characteristic function is

$$\Phi_t(\lambda) = E e^{i\lambda L_t} = e^{t\eta(\lambda)} \quad for\ each\ t \geq 0, \lambda \in \mathbb{R},$$

where
$$\eta(\lambda) = ib\lambda - a\frac{\lambda^2}{2} + \int_{\mathbb{R}\setminus\{0\}} (e^{i\lambda y} - 1 - i\lambda y I_{\{|y|<1\}})\nu(\mathrm{d}y) \qquad (6)$$

for some $b \in \mathbb{R}$, $a > 0$ and a non-negative Borel measure $\nu$ on $\mathbb{R} \setminus \{0\}$ for which
$$\int_{\mathbb{R}\setminus\{0\}} (y^2 \wedge 1)\, \nu(\mathrm{d}y) < \infty, \qquad (7)$$

or equivalently
$$\int_{\mathbb{R}\setminus\{0\}} \frac{y^2}{1+y^2}\, \nu(\mathrm{d}y) < \infty. \qquad (8)$$

Here $I_S$ is the indicator function of the set $S$. The formula (6) and the measure $\nu$ are called the characteristic exponent and the Lévy jump measure of the Lévy motion $L(t)$, respectively.

Conversely, given a mapping of the form (6), we can always construct a Lévy motion for which $\Phi_t(\lambda) = e^{t\eta(\lambda)}$ and call it a Lévy motion with the generating triplet $(b, a, \nu)$.

We reorganize the formula (6) in the form [15]:

$$\eta(\lambda) = \eta^{(1)} + \eta^{(2)} + \eta^{(3)}, \tag{9}$$

where $\eta^{(1)} = ib\lambda - a\frac{\lambda^2}{2}$, $\eta^{(2)} = \nu(\mathbb{R} \setminus (-1,1)) \int_{|y|\geq 1}(e^{i\lambda y} - 1)\frac{\nu(dy)}{\nu(\mathbb{R}\setminus(-1,1))}$, $\eta^{(3)} = \int_{0<|y|<1}(e^{i\lambda y}-1-i\lambda y)\nu(\mathrm{d}y)$. From the point of view of the characteristic function, we find that $\eta^{(1)}$ and $\eta^{(2)}$ correspond, respectively, to a linear Brownian motion with drift:

$$L_t^{(1)} = -\sqrt{a}B_t + bt, t \geq 0,$$

and a compound Poisson process:

$$L_t^{(2)} = \sum_{i=1}^{N_t} \xi_i, t \geq 0,$$

where $\{N_t : t \geq 0\}$ is a Poisson process with rate $\nu(\mathbb{R} \setminus (-1,1))$ and $\{\xi_i : i \geq 1\}$ are independent and identically distributed with distribution $\nu(dy)/\nu(\mathbb{R} \setminus (-1,1))$ concentrated on $\{y : |y| \geq 1\}$ (unless $\nu(\mathbb{R}\setminus(-1,1)) = 0$ in which case $L_t^{(2)}$ is the process which is identically zero).

The following decomposition can help us understand the Lévy motion better [15].

**The Lévy-Itô decomposition.** Given any $b \in \mathbb{R}, a \geq 0$ and measure $\nu$ concentrated on $\mathbb{R} \setminus \{0\}$ satisfying (7) or (8), there exists a probability space on which three independent Lévy processes exist, $L_t^{(1)}, L_t^{(2)}$ and $L_t^{(3)}$ where $L_t^{(1)}$ is a linear Brownian motion with drift with characteristic exponent $\eta^{(1)}$, $L_t^{(2)}$ is a compound Poisson process with characteristic exponent $\eta^{(2)}$ and $L_t^{(3)}$ is a square integrable martingale with an almost surely countable number of jumps on each finite time interval which are of magnitude less than unity and with characteristic exponent given by $\eta^{(3)}$. By taking

$$L_t = L_t^{(1)} + L_t^{(2)} + L_t^{(3)}, \tag{10}$$

we obtain that there exists a probability space on which a Lévy process is defined with characteristic exponent

$$\eta(\lambda) = ib\lambda - a\frac{\lambda^2}{2} + \int_{\mathbb{R}\setminus\{0\}}(e^{i\lambda y} - 1 - i\lambda y I_{\{|y|<1\}})\nu(\mathrm{d}y)$$

for $\lambda \in \mathbb{R}$.

The generator A of the process $L_t$ is the same as the infinitesimal generator since the Lévy process has independent and stationary increments. So the generator $A$ of the process $L_t$ is defined as $A\varphi = \lim_{t\downarrow 0}\frac{P_t\varphi-\varphi}{t}$, where $P_t \varphi(x) = E_x \varphi(L_t)$ and $\varphi$ is any function belonging to the domain of the operator $A$. Recall that the space $C_b^2(\mathbb{R})$ of $C^2$ functions with bounded derivatives up to order 2 is contained in the domain of $A$, and that for every $\varphi \in C_b^2(\mathbb{R})$

$$A\varphi(x) = b\varphi'(x) + \frac{a}{2}\varphi''(x) + \int_{\mathbb{R}\setminus\{0\}} [\varphi(x+y) - \varphi(x) - I_{\{|y|<1\}} y\varphi'(x)] \, \nu(\mathrm{d}y). \tag{11}$$

(For more details see [2, 3, 7, 24]).

For $\alpha \in (0, 2]$, a scalar symmetric $\alpha$-stable Lévy motion $L_t^\alpha$ is a Lévy process with characteristic exponent

$$\eta(\lambda) = -|\lambda|^\alpha, \lambda \in \mathbb{R}.$$

A symmetric 2-stable process is simply a Brownian motion. When $\alpha \in (0, 2)$, the generating triplet of the symmetric $\alpha$-stable process $L_t$ is $(0, 0, \nu_\alpha)$, where the Lévy jump measure $\nu_\alpha$ is given by

$$\nu_\alpha(\mathrm{d}y) = C_\alpha |y|^{-(1+\alpha)} \, \mathrm{d}y$$

with (see [7]) $C_\alpha = \dfrac{\alpha}{2^{1-\alpha}\sqrt{\pi}} \dfrac{\Gamma(\frac{1+\alpha}{2})}{\Gamma(1-\frac{\alpha}{2})}$. Sometimes, $\alpha$ is called the stability index or non-Gaussianity index.

## 4 Mean exit time, escape probability and numerical algorithms

In this section, we quantify the dynamic progression of the stochastic differential equation (5) by means of the mean exit time and the escape probability.

We first give the definition of the first exit time from the bounded domain $D = (x_1, x_3)$ as follows: the first exit time is the time needed by the density of the tumor cells $X(t)$ crossing the interval $D$, starting from an initial density $x$ in $D$, for the first time, which can be expressed by the following formula:

$$\tau_D(\omega) = \inf\{t \geq 0, X_t(\omega, x) \notin D\},$$

and the mean exit time can be denoted as $u(x) = E\tau_D(\omega)$ [5, 6, 11].

According to the theory of the infinitesimal generator [3, 19, 23, 25], we obtain that $u(x)$ satisfies the following differential-integral equation:

$$\begin{aligned} Au(x) &= -1, \quad x \in D, \\ u &= 0, \quad x \in D^c, \end{aligned} \tag{12}$$

where the generator $A$ is

$$\begin{aligned} Au &= f(x)u'(x) + \frac{a}{2}u''(x) \\ &+ \varepsilon \int_{\mathbb{R}\setminus\{0\}} [u(x+y) - u(x) - I_{\{|y|<1\}} yu'(x)] \, \nu_\alpha(\mathrm{d}y), \end{aligned} \quad (13)$$

and $D^c = \mathbb{R} \setminus D$ is the complement set of $D$.

We also consider the escape probability of the particle whose motion is described by the Eq. (5). The likelihood that the density of tumor cells, starting at a point $x$, first escapes the domain $D$ and lands in the subset $E$ of $D^c$ is called the escape probability. This escape probability, denoted by $p(x)$, satisfies the following equation [18, 22] :

$$\begin{aligned} A p(x) &= 0, \quad x \in D, \\ p|_{x \in E} &= 1, \quad p|_{x \in D^c \setminus E} = 0, \end{aligned} \quad (14)$$

where $A$ is the generator defined in (13). In this paper, $E = (-\infty, x_1]$ means the probability of tumor extinction, while $E = [x_3, +\infty)$ means the probability of becoming malignant.

As is seen above, the non-Gaussianity of the Lévy noise results in the nonlocality of the equations of the mean exit time and the escape probability.

Now we describe the numerical algorithms of the equation (12) for the non-symmetric interval (base on the scheme in the reference [11], which has given the algorithms for the symmetric interval case). To make it apply to the general case, we use $D = (c, d), c \leq 0, d > 0$ in the process of discretization instead of $D = (x_1, x_3)$, so the equation (12) can be recast as:

$$\frac{a}{2}u''(x) + f(x)u'(x) + \varepsilon C_\alpha \int_{\mathbb{R}\setminus\{0\}} \frac{u(x+y) - u(x) - I_{\{|y|<\delta\}} yu'(x)}{|y|^{1+\alpha}} \, \mathrm{d}y = -1, \quad (15)$$

for $x \in (c, d)$; and $u(x) = 0$ for $x \notin (c, d)$.

In order to make the integrand $u(x + y)$ vanish, we divide the integral interval as $\int_\mathbb{R} = \int_{-\infty}^{c-x} + \int_{c-x}^{d-x} + \int_{d-x}^{+\infty}$ and choose $\delta = \min\{|c - x|, |d - x|\}$, then we obtain

$$\begin{aligned} &\frac{a}{2}u''(x) + f(x)u'(x) - \frac{\varepsilon C_\alpha}{\alpha}\left[\frac{1}{(x-c)^\alpha} + \frac{1}{(d-x)^\alpha}\right]u(x) \\ &+ \varepsilon C_\alpha \int_{c-x}^{d-x} \frac{u(x+y) - u(x) - I_{\{|y|<\delta\}} yu'(x)}{|y|^{1+\alpha}} \, \mathrm{d}y = -1, \end{aligned} \quad (16)$$

for $x \in (c, d)$; and $u(x) = 0$ for $x \notin (c, d)$.

According to the property of the indicator function, the Eq. (16) can be rewritten as:

$$\frac{a}{2}u''(x) + f(x)u'(x) - \frac{\varepsilon C_\alpha}{\alpha}\left[\frac{1}{(x-c)^\alpha} + \frac{1}{(d-x)^\alpha}\right]u(x) \qquad (17)$$

$$+\varepsilon C_\alpha \int_{c-x}^{x-d} \frac{u(x+y) - u(x)}{|y|^{1+\alpha}}\,\mathrm{d}y + \varepsilon C_\alpha \int_{x-d}^{d-x} \frac{u(x+y) - u(x) - yu'(x)}{|y|^{1+\alpha}}\,\mathrm{d}y = -1,$$

for $x \geq \frac{c+d}{2}$, and

$$\frac{a}{2}u''(x) + f(x)u'(x) - \frac{\varepsilon C_\alpha}{\alpha}\left[\frac{1}{(x-c)^\alpha} + \frac{1}{(d-x)^\alpha}\right]u(x) \qquad (18)$$

$$+\varepsilon C_\alpha \int_{x-c}^{d-x} \frac{u(x+y) - u(x)}{|y|^{1+\alpha}}\,\mathrm{d}y + \varepsilon C_\alpha \int_{c-x}^{x-c} \frac{u(x+y) - u(x) - yu'(x)}{|y|^{1+\alpha}}\,\mathrm{d}y = -1,$$

for $x < \frac{c+d}{2}$.

Now, let us take the appropriate stepsize $h$, so that $\frac{c}{h}, \frac{d}{h}, \frac{c+d}{2h}$ are integers, and define $x_j = jh$ for $\frac{c-d}{h} \leq j \leq \frac{d-c}{h}$. We use the notation $U_j$ to indicate the numerical solution of $u$ at $x_j$. Then, the differential-integral equations (17) and (18) can be discretized, respectively, using the central difference scheme for derivatives and "punched-hole" trapezoidal rule

$$\frac{a}{2}\frac{U_{j-1} - 2U_j + U_{j+1}}{h^2} + f(x_j)\frac{U_{j+1} - U_{j-1}}{2h} - \frac{\varepsilon C_\alpha}{\alpha}\left[\frac{1}{(x_j - c)^\alpha} + \frac{1}{(d - x_j)^\alpha}\right]U_j$$

$$+ \varepsilon C_\alpha h \sum_{k=\frac{c}{h}-j}^{j-\frac{d}{h}}{}'' \frac{U_{j+k} - U_j}{|x_k|^{1+\alpha}} + \varepsilon C_\alpha h \sum_{k=j-\frac{d}{h}, k\neq 0}^{\frac{d}{h}-j}{}'' \frac{U_{j+k} - U_j - (U_{j+1} - U_{j-1})x_k/2h}{|x_k|^{1+\alpha}} = -1,$$

$$(19)$$

where $j = \frac{c+d}{2h}, \frac{c+d}{2h} + 1, \cdots, \frac{d}{h} - 1$. The meaning of the modified summation symbol $\sum''$ is that the quantities corresponding to the two endpoints of the integral interval should be multiplied by $1/2$.

$$\frac{a}{2}\frac{U_{j-1} - 2U_j + U_{j+1}}{h^2} + f(x_j)\frac{U_{j+1} - U_{j-1}}{2h} - \frac{\varepsilon C_\alpha}{\alpha}\left[\frac{1}{(x_j - c)^\alpha} + \frac{1}{(d - x_j)^\alpha}\right]U_j$$

$$+ \varepsilon C_\alpha h \sum_{k=j-\frac{c}{h}}^{\frac{d}{h}-j}{}'' \frac{U_{j+k} - U_j}{|x_k|^{1+\alpha}} + \varepsilon C_\alpha h \sum_{k=\frac{c}{h}-j, k\neq 0}^{j-\frac{c}{h}}{}'' \frac{U_{j+k} - U_j - (U_{j+1} - U_{j-1})x_k/2h}{|x_k|^{1+\alpha}} = -1,$$

$$(20)$$

where $j = \frac{c}{h} + 1, \frac{c}{h} + 2, \cdots, \frac{c+d}{2h} - 1$. The boundary conditions require that the values of $U_j$ vanish when the index $j \leq \frac{c}{h}$ or $j \geq \frac{d}{h}$.

The truncation errors of the central difference scheme for derivatives in (19) and (20) are of 2nd-order $O(h^2)$. According to the error analysis of the reference [11], the leading-order error of the quadrature rule is $-\zeta(\alpha - 1)u''(x)h^{2-\alpha} + O(h^2)$, where $\zeta$ is the Riemann zeta function. Thus, the modified results with 2nd-order accuracy can be expressed as follows:

$$C_h \frac{U_{j-1} - 2U_j + U_{j+1}}{h^2} + f(x_j)\frac{U_{j+1} - U_{j-1}}{2h} - \frac{\varepsilon C_\alpha}{\alpha}\left[\frac{1}{(x_j - c)^\alpha} + \frac{1}{(d - x_j)^\alpha}\right] U_j$$

$$+ \varepsilon C_\alpha h \sum_{k=\frac{c}{h}-j}^{j-\frac{d}{h}}{}'' \frac{U_{j+k} - U_j}{|x_k|^{1+\alpha}} + \varepsilon C_\alpha h \sum_{k=j-\frac{d}{h},k\neq 0}^{\frac{d}{h}-j}{}'' \frac{U_{j+k} - U_j - (U_{j+1} - U_{j-1})x_k/2h}{|x_k|^{1+\alpha}} = -1,$$
(21)

for any $0 < \alpha < 2$, $j = \frac{c+d}{2h}, \frac{c+d}{2h}+1, \cdots, \frac{d}{h}-1$. Here $C_h = \frac{a}{2} - \varepsilon C_\alpha \zeta(\alpha - 1)h^{2-\alpha}$.
Similarly, for $j = \frac{c}{h} + 1, \frac{c}{h} + 2, \cdots, \frac{c+d}{2h} - 1$,

$$C_h \frac{U_{j-1} - 2U_j + U_{j+1}}{h^2} + f(x_j)\frac{U_{j+1} - U_{j-1}}{2h} - \frac{\varepsilon C_\alpha}{\alpha}\left[\frac{1}{(x_j - c)^\alpha} + \frac{1}{(d - x_j)^\alpha}\right] U_j$$

$$+ \varepsilon C_\alpha h \sum_{k=j-\frac{c}{h}}^{\frac{d}{h}-j}{}'' \frac{U_{j+k} - U_j}{|x_k|^{1+\alpha}} + \varepsilon C_\alpha h \sum_{k=\frac{c}{h}-j,k\neq 0}^{j-\frac{c}{h}}{}'' \frac{U_{j+k} - U_j - (U_{j+1} - U_{j-1})x_k/2h}{|x_k|^{1+\alpha}} = -1.$$
(22)

$U_j = 0$ for $j \leq \frac{c}{h}$ or $j \geq \frac{d}{h}$.

Thus, we can get the mean exit time by solving the linear equations (21) and (22).

When it mentions to the escape probability, the Eq. (14) can be written as follows:

$$\frac{a}{2}p''(x) + f(x)p'(x) - \frac{\varepsilon C_\alpha}{\alpha}\left[\frac{1}{(x - c)^\alpha} + \frac{1}{(d - x)^\alpha}\right] p(x)$$

$$+\varepsilon C_\alpha \int_{c-x}^{d-x} \frac{p(x+y) - p(x) - I_{\{|y|<\delta\}}yp'(x)}{|y|^{1+\alpha}} \, dy = -\frac{\varepsilon C_\alpha}{\alpha}\left[\frac{1}{(x-c)^\alpha}\right], \quad (23)$$

for $x \in D = (c, d)$. The boundary conditions are $p(x) = 0$ for $x \in [d, +\infty)$, and $p(x) = 1$ for $x \in (-\infty, c]$ in this paper. By the similar process above, we can obtain the escape probability.

# 5 Stochastic evolution of tumor density

We fix the parameters $\theta = 0.1, \beta = 3.0$, as in [10], and focus on the impact of the noise parameters on the mean exit time and the escape probability.

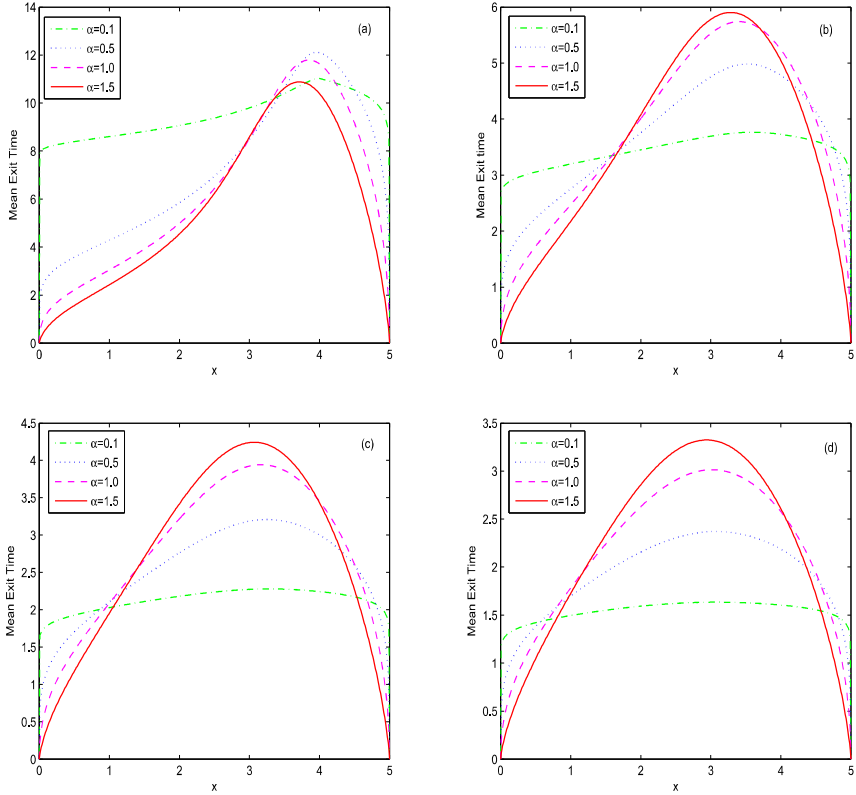

Figure 2: Mean exit time $u(x)$: (a) $a = 0.0, \varepsilon = 0.1$. (b) $a = 0.0, \varepsilon = 0.3$. (c) $a = 0.0, \varepsilon = 0.5$. (d) $a = 0.0, \varepsilon = 0.7$.

With these parameters, deterministic stable states are $x_1 = 0$ and $x_3 = 5$. The interval $D = (0, 5)$ encloses the tumor density from 0 (state of tumor extinction) to 5 (state of stable tumor).

## 5.1 Mean exit time

In Figures 2-6, we plot the mean exit time of the tumor density $X_t$ from the interval $D = (0, 5)$, between the state of tumor extinction and the state of stable tumor. For example, $u(2.5)$ is the mean time that the system remains in $D$ starting from tumor density $x = 2.5$, before 'exiting' to outside $D$. All figures are color-coded, as well as coded with distinct line patterns, for various parameters. They have better online visibility.

In Figures 2 and 4, the diffusion coefficient $a = 0$ means the noise is a pure non-Gaussian noise, while $a > 0$ means the Gaussian noise is combined with non-Gaussian noise. They demonstrate the dependence of the mean exit time on the stability index $\alpha$, when the diffusion coefficient $a$ and the noise intensity $\varepsilon$ are fixed. In fact, when the noise intensity $\varepsilon$ is very small ($\varepsilon < 0.1$), the mean exit time decreases as the stability index $\alpha$ becomes large. As $\varepsilon$ gradually increases and is fixed, in the different stages of tumor, the relationships between the mean exit time and the stability index $\alpha$ are different. In the early tumor stage (i.e., small tumor density $x$), the mean exit time decreases with increasing stability index $\alpha$. When the tumor density $x$ increases to a certain extent, this relationship becomes opposite. But if the density of the tumor cells continues to increase, the mean exit time decreases with increasing stability index $\alpha$ again. Moreover, with the increasing noise intensity $\varepsilon$, the mean exit time at the middle tumor stage becomes bigger.

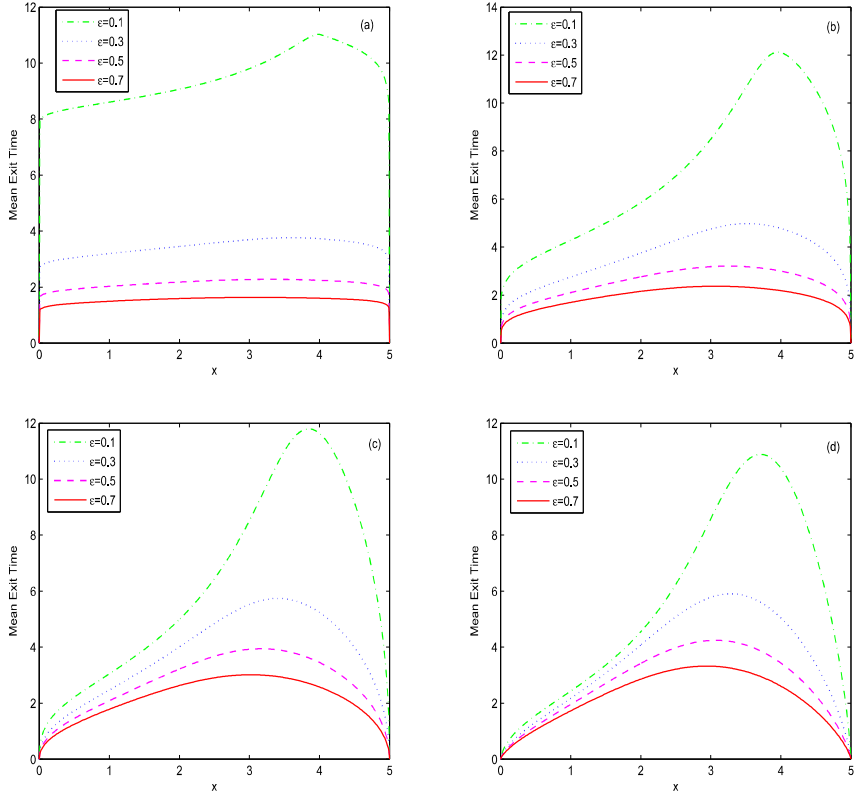

Figure 3: Mean exit time $u(x)$: (a) $a = 0.0, \alpha = 0.1$. (b) $a = 0.0, \alpha = 0.5$. (c) $a = 0.0, \alpha = 1.0$. (d) $a = 0.0, \alpha = 1.5$.

When the diffusion coefficient $a$ and stability index $\alpha$ are fixed, the mean exit time becomes shorter with the increase of the noise intensity $\varepsilon$ for all the different initial tumor densities (see Figures 3 and 5). Besides, when the stability index $\alpha$ is small, raising the noise intensity $\varepsilon$ makes the mean exit time becomes shorter and even roughly constant (see Figure 3(a)), which means for all the tumor densities, the lifetimes staying in the domain $D$ are almost the same.

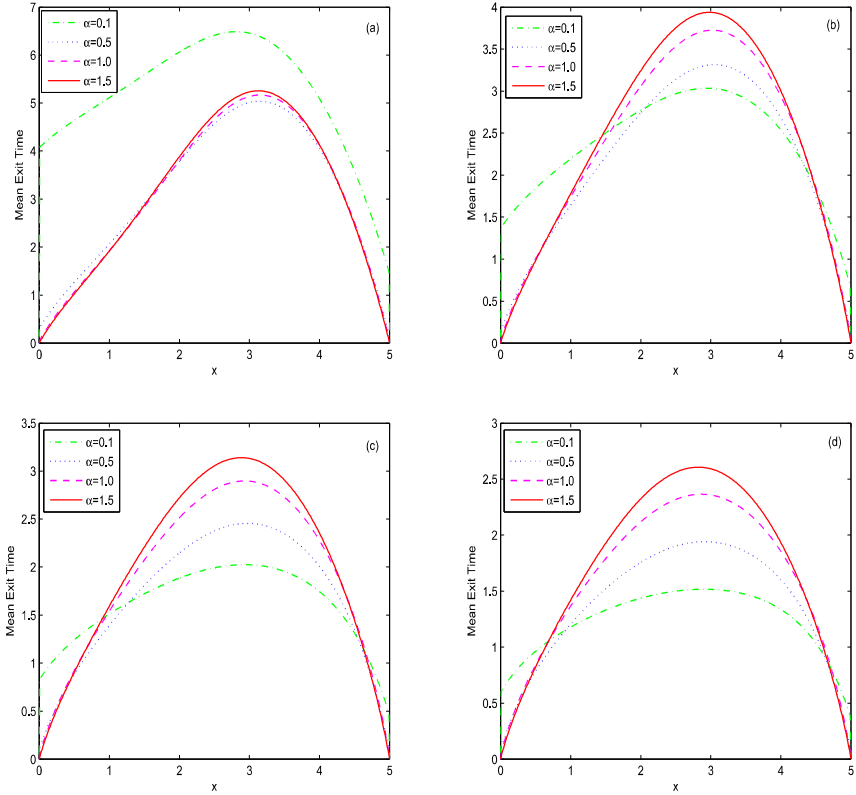

Figure 4: Mean exit time $u(x)$:(a) $a = 0.5, \varepsilon = 0.1$. (b) $a = 0.5, \varepsilon = 0.3$.(c) $a = 0.5, \varepsilon = 0.5$. (d) $a = 0.5, \varepsilon = 0.7$.

In Figure 6, we find that when the noise parameters $\alpha$ and $\varepsilon$ are fixed, in the presence of the Gaussian noise (non-zero $a$ values), the lifetimes of the tumor with all densities become shorter and decrease as the Gaussian diffusion coefficient $a$ increases.

## 5.2 Escape probability

The mean exit time $u(x)$ quantifies the mean lifetime for a tumor, with density value $x$, staying in the interval $D = (0, 5)$. After the mean exit time, the system exits $D$, either to the left or to the right. But this exiting direction is not distinguished in the mean exit time. We would also like to understand the tumor evolution after this stage, i.e., will it reach the tumor-free state (exit to the left of $D = (0, 5)$) or the malignant state (exit to the right of $D = (0, 5)$)? Note that $x = 0$ means tumor density is 0, while $x = 5$ indicates the tumor density is the highest, for the biomedical parameters we have chosen.

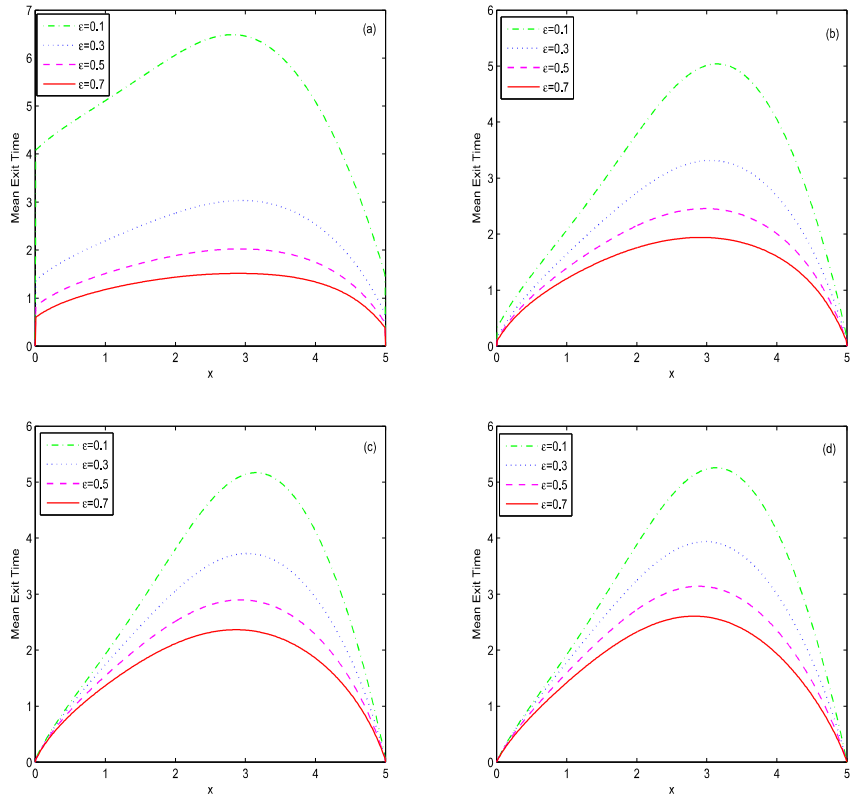

Figure 5: Mean exit time $u(x)$: (a) $a = 0.5, \alpha = 0.1$. (b) $a = 0.5, \alpha = 0.5$. (c) $a = 0.5, \alpha = 1.0$. (d) $a = 0.5, \alpha = 1.5$.

Biologically speaking, we are more interested in the likelihood of the tumor extinction, induced by the noise after the mean exit time. To this end, we now compute the escape probability $p(x)$ for the system to escape $D = (0, 5)$ to the left, showed in Figures 7-11. For example, $p(2.5)$ is the likelihood that the system at tumor density $x = 2.5$ becomes tumor-free (i.e., tumor cells become extinct), under a given set of biomedical and noise parameters.

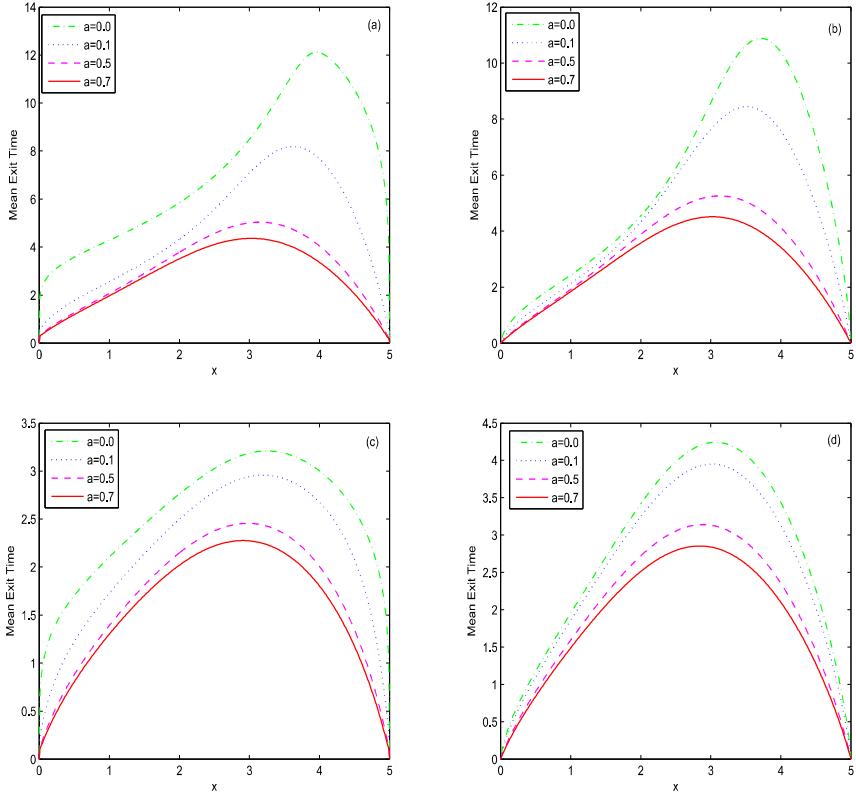

Figure 6: Mean exit time $u(x)$: (a) $\varepsilon = 0.1, \alpha = 0.5$. (b) $\varepsilon = 0.1, \alpha = 1.5$. (c) $\varepsilon = 0.5, \alpha = 0.5$. (d) $\varepsilon = 0.5, \alpha = 1.5$.

*Pure non-Gaussian noise case.*
First, we consider the pure jump noise case, i.e., the diffusion coefficient $a = 0$.

Figure 7 shows the dependence of the escape probability $p(x)$ on the stability index $\alpha$ under different noise intensities $\varepsilon$. We find that there exists a critical tumor density value $I_\varepsilon$, depending on noise intensity $\varepsilon$. When the density of tumor cells is below $I_\varepsilon$, the escape probability increases as the stability index $\alpha$ increases. That means that the bigger the stability index $\alpha$ is, the higher the probability of the tumor cells becoming extinct. While when the density of tumor cells is above $I_\varepsilon$, the escape probability decreases as the stability index $\alpha$ increases. So for larger stability index $\alpha$, the tumor density is more likely to grow. Moreover, the critical value $I_\varepsilon$ decreases with $\varepsilon$.

Figure 8 shows how the escape probability $p(x)$ varies with the noise intensity $\varepsilon$, for a few different stability indexes $\alpha$. We observe that the escape probability decreases as the noise intensity $\varepsilon$ increases, in general. This appears to indicate that large jump noise intensity has an adverse impact on tumor development.

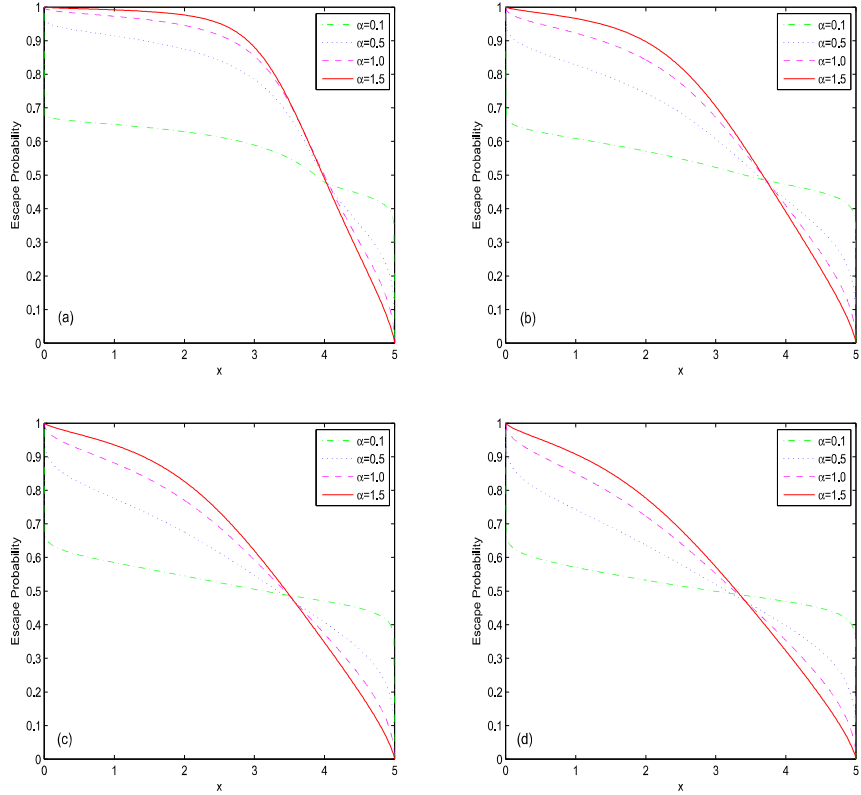

Figure 7: Escape probability $p(x)$: (a) $a = 0.0, \varepsilon = 0.1$. (b) $a = 0.0, \varepsilon = 0.3$. (c) $a = 0.0, \varepsilon = 0.5$. (d) $a = 0.0, \varepsilon = 0.7$.

*Combined non-Gaussian and Gaussian noise case.*
Now, we discuss the effect of combined non-Gaussian and Gaussian noise on the escape probability through the left boundary $x = 0$. In this case, both diffusion coefficient $a$ and jump coefficient $\varepsilon$ are non-zero. We compute the likelihood that tumor density becomes zero (cancer extinction).

Comparing Figures 7 with 9, and 8, with 10, respectively, we find that the additional Gaussian noise does not change much of the overall evolution tendency of the escape probability, although the critical value $I_\varepsilon$ has changed a little bit in Figure 10.

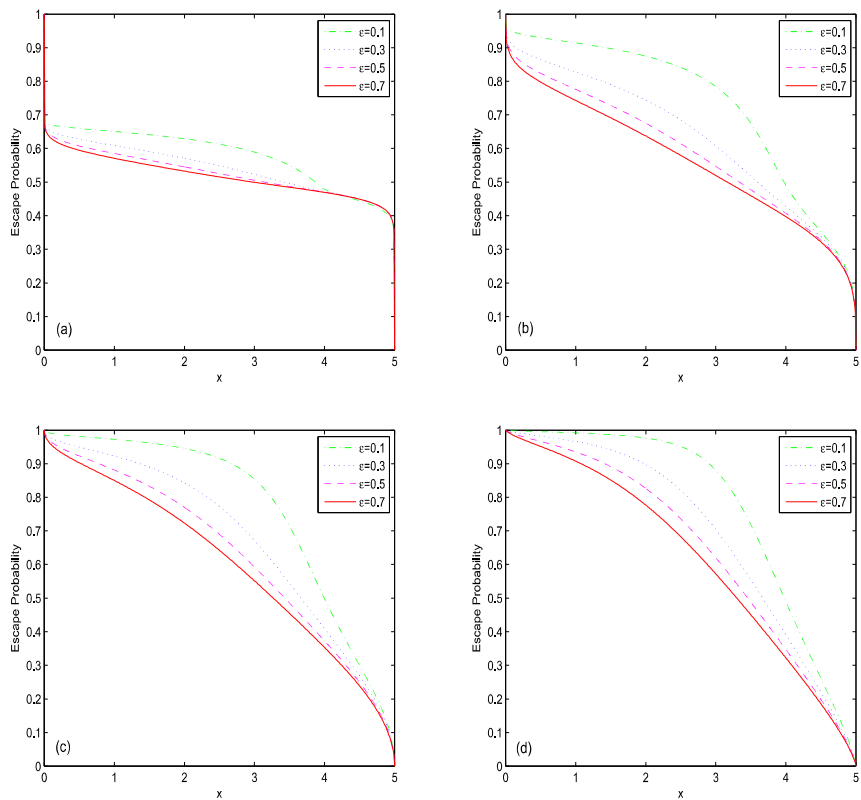

Figure 8: Escape probability $p(x)$: (a) $a = 0.0, \alpha = 0.1$. (b) $a = 0.0, \alpha = 0.5$. (c) $a = 0.0, \alpha = 1.0$. (d) $a = 0.0, \alpha = 1.5$.

As shown in Figure 11, when the stability index $\alpha > 1$, the escape probability is reduced due to increasing Gaussian noise (i.e., increasing diffusion coefficient $a$). For the stability index $\alpha < 1$, there is a critical value $I_a$, depending on diffusion coefficient $a$, so that the likelihood that the tumor becomes extinct increases with $a$ for initial tumor density $x < I_a$, but decreases with $a$ for initial tumor density $x > I_a$.

# 6 Conclusion

We have studied the dynamical evolution of a tumor growth model with immunization, under non-Gaussian random influences. We focus on quantifying the mean lifetime for the tumor density remains between the tumor-free state and the state of stable tumor, and computing the likelihood that a tumor with a certain initial density becomes extinct (i.e., tumor-free), under a set of medical and noise parameters.

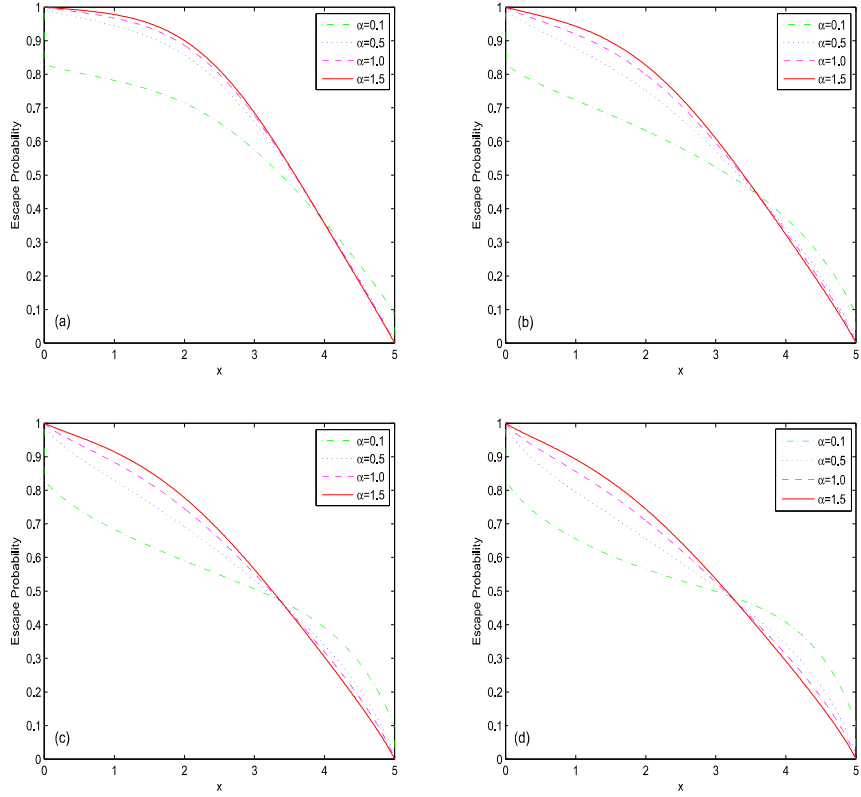

Figure 9: Escape probability $p(x)$:(a) $a = 0.5, \varepsilon = 0.1$. (b) $a = 0.5, \varepsilon = 0.3$.(c) $a = 0.5, \varepsilon = 0.5$. (d) $a = 0.5, \varepsilon = 0.7$.

From our numerical experiments, it is observed that the mean lifetime and the escape probability depend strongly on the initial tumor density and the parameters of non-Gaussian Lévy noise. For the inchoate patients, the smaller jumps with higher jump probability (i.e., $\alpha$ is closer to 2) facilitate the control of the tumor cells. But in the late tumor stage, the larger jumps with lower frequency (that is to say $0 < \alpha < 1$) are better to slow the cancerous progression. Besides, although the larger amount of the noise can shorten the mean lifetime, it reduces the probability that induces the tumor extinction at the same time. Due to the complexity of the interplay between the tumor issue and the immune cells, any small fluctuations from the environment may influence the progression of the tumor cells considerably.

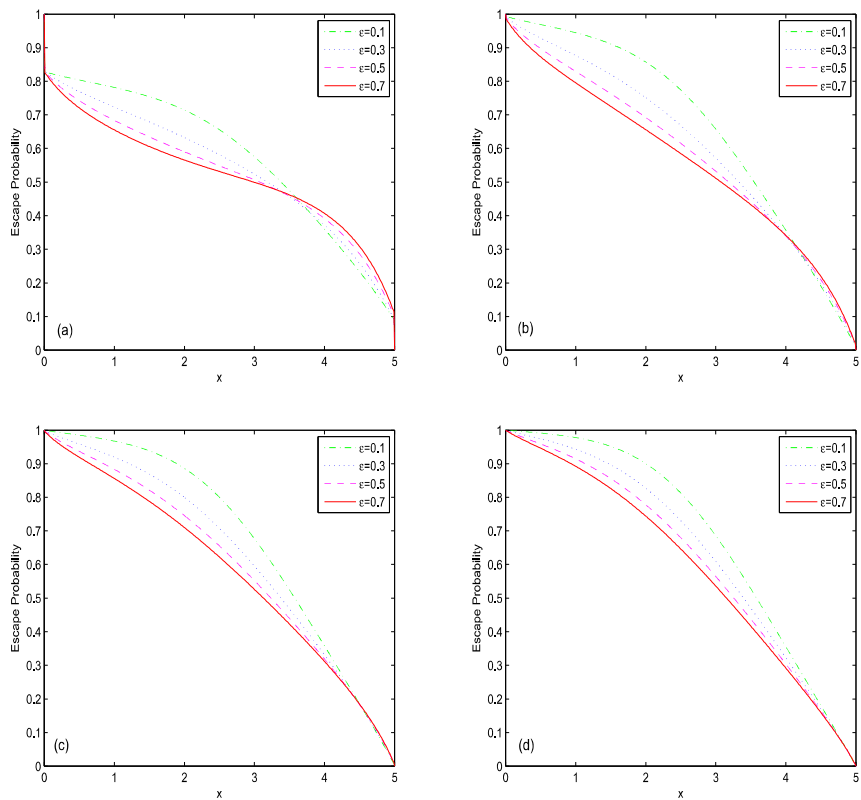

Figure 10: Escape probability $p(x)$: (a) $a = 0.5, \alpha = 0.1$. (b) $a = 0.5, \alpha = 0.5$. (c) $a = 0.5, \alpha = 1.0$. (d) $a = 0.5, \alpha = 1.5$.

**Acknowledgements.** We thank Jian Ren and Caishi Wang for helpful discussions.

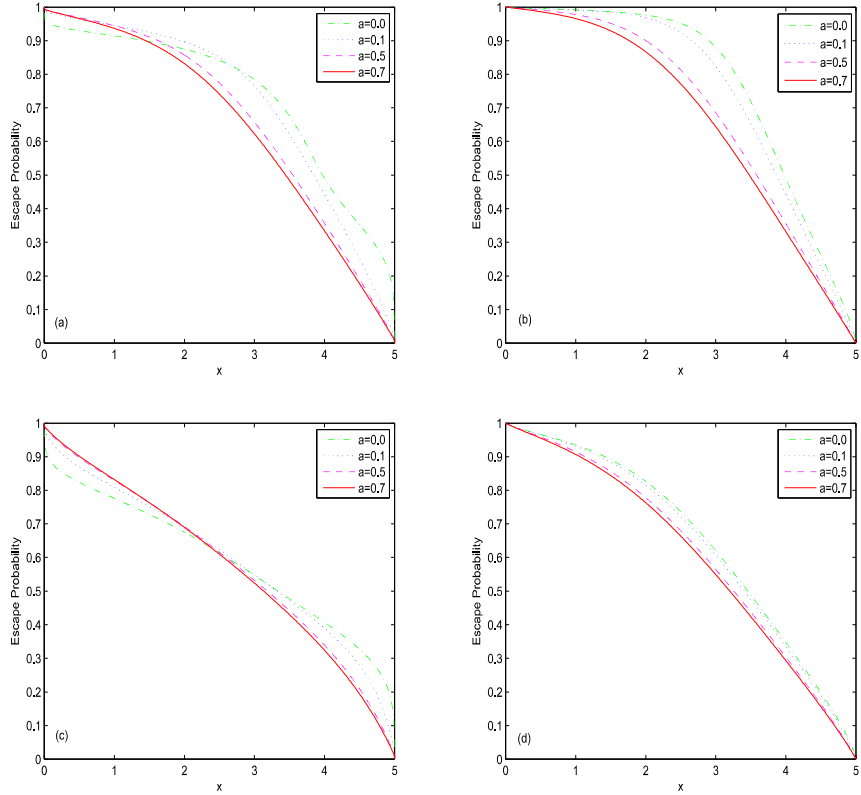

Figure 11: Escape probability $p(x)$: (a) $\varepsilon = 0.1, \alpha = 0.5$. (b) $\varepsilon = 0.1, \alpha = 1.5$. (c) $\varepsilon = 0.5, \alpha = 0.5$. (d) $\varepsilon = 0.5, \alpha = 1.5$.